\documentclass[12pt]{article}
\usepackage{amsfonts}
\usepackage{amssymb}
\usepackage{mathrsfs}
\usepackage{float}
\usepackage{fancyhdr}
\usepackage{amsmath,,amsthm}
\usepackage{amsfonts,amssymb}
\usepackage{graphicx,color}
\usepackage{booktabs,tabularx}
\usepackage{dsfont}
\usepackage{epsfig}
\usepackage{graphicx}
\usepackage{subfigure}

\date{ }

\title{The Periodicity to a Kind of Generalized Collatz Problem}
\author{Sensen Chen, Qing-You Sun\footnote{Corresponding author:
qysun@hznu.edu.cn} \hspace{0.2cm}and Yushu Zhu
\\ \small Hangzhou Normal University, Hangzhou
311121, China }

\begin{document}
\maketitle
\begin{abstract}
The Collatz problem is related to the fixed point problem, and is
widely used in mathematics. It has attracted a wide range of math
enthusiasts, but is still difficult to solve. So, this article aimed
to study the extension of the Collatz problem, more widely, in
$k$-adic. We define a new sequence called $\mathcal{Z}$
transformation sequence. Under a suitable assumptions, we can prove
that the limit set of the $\mathcal{Z}$ transformation sequence must
be $M=\{1,2\}$.

\vskip 3mm

\noindent\textbf{Key words and phrases}: Collatz problem; $k$-adic;
$\mathcal{Z}$ transformation; limit set

\noindent\textbf{Mathematics Subject Classification 2010}: 11A99,
11A67, 11B83

\end{abstract}

\baselineskip=5mm

\section{Introduction}

\hspace{0.5cm} The Collatz problem has been widely studied in the
past 100 years, and many achievements with great value have been
obtained, although the Collatz problem cannot be effectively solved.
However, this does not affect the important position of the Collatz
problem in mathematics.

The essence of the Collatz problem is a fixed-point problem, and its
research results are important for solving fixed-point problems.
Fixed-point problems are an important foundation for modern analysis
and topology, so the study of the Collatz problem is of great value.
Moreover, it is widely used and has important connection with
dynamic system, fractal geometry and other fields. It is also the
theoretical basis of cryptography research.

The Collatz problem was proposed by Collatz \cite {1}, a German
mathematician, in 1937. He conjectured the following
number-theoretic function

\begin{equation} \label{eq:1}
f(x)=\left\{\begin{array}{c}{\displaystyle\frac{x}{2}, \text { when
} x \text { is even, }}
\\ {\displaystyle\frac{3 x + 1}{2}, \text { when } x \text { is odd.
}}\end{array}\right.
\end{equation}

There exists a finite positive integer $k$, s.t.,  $f^{n}
(x)\in\{1,4,2\}$ when $n\geq k$, and $x\in\mathds{N}$. We call this
problem Collatz conjecture. This problem mainly studies the
periodicity of sequence transformation.

The Collatz problem can be stated in the form of $2$-adic \cite {2}.
The integer in $\mathds{Z}_{2}$ can be expressed as

\begin{equation*}
\alpha=a_{m-1} 2^{m-1}+a_{m-2} 2^{m-2}+ \cdots +a_{1} \cdot 2 +
a_{0},
\end{equation*}
where $a_{i}=0$, or $1$, $0\leq i<m-1,1\leq m<\infty$. It called the
$2$-adic form of $\alpha$.

One can define congruence on $\mathds{Z}_{2}$ by $\alpha \equiv
\beta\left(\bmod 2^{m}\right)$ if the first $m$ $2$-adic digits of
$\alpha$ and $\beta$ agree. Addition and multiplication on
$\mathds{Z}_{2}$ are given by

\begin{equation*}
\begin{array}{c}{(X)_{2}=(\alpha+\beta)_{2}=X\left(\bmod 2^{k}\right)
=\alpha\left(\bmod 2^{k}\right)+\beta\left(\bmod 2^{k}\right)},
\\ {\quad(X)_{2}=(\alpha \beta)_{2}=X\left(\bmod 2^{k}\right)
=\alpha\left(\bmod 2^{k}\right) \cdot \beta\left(\bmod
2^{k}\right).}\end{array}
\end{equation*}

Now, one can extend the definition of the function
$T:\mathds{Z}\rightarrow\mathds{Z}$ given by \eqref{eq:1} to
$T:\mathds{Z}_{2}\rightarrow\mathds{Z}_{2}$ by

\begin{equation} \label{eq:2}
T(x)=\left\{\begin{array}{c}{\displaystyle\frac{x}{2}, x \equiv
0(\bmod 2)},
\\ {\displaystyle\frac{3 x+1}{2},
x \equiv 1(\bmod 2)}.\end{array}\right.
\end{equation}

E. Heppner \cite{3} got an important corollary and extended it to
the form of module $p$. In \cite{3}, the following mapping is given

\begin{equation*}
T(n)=T_{m, p, j}(n) \left\{\begin{array}{c}{\displaystyle\frac{m
n+r_{j}}{p}, n \equiv j(\bmod p)},
\\ {\displaystyle\frac{n}{p}, n
\equiv 0(\bmod p)},\end{array}\right.
\end{equation*}

where $(m,p)=1$ with $p>2$, and $r_j \equiv -m_j (\bmod p)$.
Meanwhile, the transform are shown to depend on relative sizes of
$m$ and $p$.

Keith R. Matthews \cite{4} studied some maps on the rings of
integers in an algebraic number field, such as

\begin{equation*}
U(\alpha)=\left\{\begin{array}{c}{\displaystyle\frac{(1-\sqrt{2})
\alpha}{\sqrt{2}}, \alpha \equiv 0(\bmod 2)},
\\ {\displaystyle\frac{3
\alpha+1}{\sqrt{2}}, \alpha \equiv 1(\bmod 2)}.\end{array}\right.
\end{equation*}

Certainly, there are many other research results of such problems
and their extended problems (\cite{5}-\cite{8}), especially through
computers to solve such problems. This article mainly explores a
more extensive periodic problem of a particular sequence
transformation based on the Collatz problem.

In the present article, we will extend the Collatz problem in
p-adic. Then, we give a new transformation, called $\mathcal{Z}$
transformation. Through the study of the property of the
$\mathcal{Z}$ transformation, we find that under some suitable
assumptions, the $\mathcal{Z}$ transformation has a periodic
characteristic.

\section{Definition and main result}

\hspace{0.5cm} For any given $m$-digit natural number $n$, it can be
expressed in $k$-adic as

\begin{equation} \label{eq:n}
n=a_{m-1} k^{m-1}+a_{m-2} k^{m-2}+\dots+a_{1} k+a_{0},
\end{equation}
where $a_i \in \mathds{Z}_0^+$, $0 \leq a_i < k$ $(0 \leq i \leq
m-1)$, and $a_{m-1} \neq 0$.

We introduce the function
\begin{equation} \label{eq:f}
f(x)=\left\{\begin{array}{cc}{\displaystyle\frac{(x+p-1)(x+2 p-1)}{p
\cdot p}} & {x \equiv 1(\bmod p)},
\\ {\displaystyle\frac{x+p-2}{p}} & {x \equiv 2(\bmod p)},
\\\vdots
\\{\displaystyle\frac{x+1}{p}} & {x = p-1(\bmod p)},
\\ {\displaystyle\frac{x}{p}} & {x \equiv
0(\bmod p)}.\end{array}\right.
\end{equation}

\noindent\textbf{Definition 1}\quad \emph{ Let
$\mathcal{Z}_k(n)=\displaystyle\sum_{i=0}^{m-1}f(a_i)$, where $f(x)$
is defined as \eqref{eq:f}. Then, $\mathcal{Z}_k(n)$ is called the
$\mathcal{Z}$ transformation of $n$ in $k$-adic.}

\noindent\textbf{Definition 2}\quad \emph{ Denote $n_j =
\mathcal{Z}_k^{j} (n), \quad (j\in \mathds{Z}_0^+)$. It is clear
that, $\{n_j\}$ is a sequence, it is called $\mathcal{Z}$
transformation sequence in $k$-adic.}

We give our main result as the following Theorem.

\noindent\textbf{Theorem 1}\quad \emph{Under the assumptions
\begin{equation} \label{eq:H}
p+2\leq k < p^2 -3p+2, \quad and \quad p>5,
\end{equation}
the $\mathcal{Z}$ transformation sequence in $k$-adic has a limit
set $M=\{1,2\}$, i.e., there exists a finite positive integer
$\lambda$, subject to $n_\mu \in\{1,2\}$ when $\mu \geq \lambda$.}

\section{The proof of Theorem 1}

\hspace{0.5cm} To prove the Theorem, we give two Lemmas first.

\noindent\textbf{Lemma 1}\quad \emph{The result of adding two
$k$-adic numbers is the same as the result of adding them in decimal
to become $k$-adic.}

This conclusion is obvious (see \cite{2}). We omit it.

\noindent\textbf{Lemma 2}\quad \emph{For any given $m$-digit natural
number $n$, which is expressed in $k$-adic as \eqref{eq:n}, under
the assumptions \eqref{eq:H}, $\mathcal{Z}_k(n)$ is a integer not
exceeding $(m-1)$-dgit in $k$-adic when $m\geq 3$.}

\begin{proof}
Denote
\begin{equation} \label{eq:r}
k = rp +s +1,
\end{equation}
where $r,s\in \mathds{N}$ and $1\leq s\leq p$.

Then, by calculating the assumptions \eqref{eq:H} are equivalent to
\begin{equation} \label{eq:H1}
1 \leq r \leq p-4,
\end{equation}
or
\begin{equation} \label{eq:H2}
(r+1)(r+2) \leq rp+1.
\end{equation}

Due to $0\leq a_i \leq k-1$, by the definition of $f(x)$ in
\eqref{eq:f}, it obtains
\begin{equation*}
f(a_i)\leq f(rp+1)=(r+1)(r+2).
\end{equation*}
Thus, $\mathcal{Z}_k(n) \leq m(r+1)(r+2)$.

For $m=3$, by \eqref{eq:H}, \eqref{eq:r}, and \eqref{eq:H2}, we have
\begin{equation*}
\mathcal{Z}_k(n) \leq 3(r+1)(r+2) < 3(rp+1)=3(k-s)<3k<k^{3-1} .
\end{equation*}
For $m=i$, we assume that it holds
\begin{equation*}
\mathcal{Z}_k(n) \leq i(r+1)(r+2) < k^{i-1} .
\end{equation*}
Then, for $m=i+1$, it's obviously that
\begin{equation*}
\mathcal{Z}_k(n) \leq (i+1)(r+1)(r+2) < k^{i-1}\cdot
\displaystyle\frac{i+1}{i} < k^{(i+1)-1} .
\end{equation*}
Therefore, by mathematical induction, we have
\begin{equation}
\mathcal{Z}_k(n) < k^{m-1} ,
\end{equation}
where $m \geq 3$. That means, $\mathcal{Z}_k(n)$ is a integer not
exceeding $(m-1)$-dgit in $k$-adic when $m\geq 3$. Thus, the proof
of Lemma 2 is done.

\end{proof}

\noindent{\bf \emph{Proof of Theorem 1.}} For any given $m$-digit
natural number $n$, which is expressed in $k$-adic as \eqref{eq:n},
we will prove the Theorem in three cases.

{\bf Case 1.} $n$ is a single-digit number.

In this case, it is easy to show that $1\leq n < k$, and
$\mathcal{Z}_k(n) = f(n)$.

If $n \not\equiv 1(\bmod p)$, by the definition of $f$, that is
\eqref{eq:f}, we can easily get $\mathcal{Z}_k(n) = f(n) < n$. It
implies that the $\mathcal{Z}$ transformation make
$\mathcal{Z}_k(n)$ smaller than $n$, and $\mathcal{Z}_k(n)$ is still
a single-digit number. Therefore, there exists a nonnegative integer
$\lambda$, which satisfies $\mathcal{Z}_k^{\lambda}(n) =1$, or
$\mathcal{Z}_k^{\lambda}(n) \equiv 1(\bmod p)$. The second case is
just what we will discuss next.

If $n \equiv 1(\bmod p)$ and $n \neq 1$, let $n=\tilde{r}p+1$. So,
by using $n < k$ and \eqref{eq:r}, we have $1\leq \tilde{r} \leq r$.\\
Noting the definition of $f$, it can obtain that
\begin{equation*}
\mathcal{Z}_k(n) = f(n) = f(\tilde{r}p+1) =
(\tilde{r}+1)(\tilde{r}+2).
\end{equation*}
Then, under the equivalent assumption \eqref{eq:H1} and
\eqref{eq:H2},
\begin{equation*}
\mathcal{Z}_k(n) < \tilde{r}p+1 = n,
\end{equation*}
here, $1\leq \tilde{r} \leq r \leq p-4$. Also, the $\mathcal{Z}$
transformation make $\mathcal{Z}_k(n)$ smaller than $n$, and
$\mathcal{Z}_k(n)$ is still a single-digit number. It means, there
exists a nonnegative integer $\lambda$, which satisfies
$\mathcal{Z}_k^{\lambda}(n) =1$, or $\mathcal{Z}_k^{\lambda}(n)
\not\equiv 1(\bmod p)$. The second case is exactly what we discussed
earlier.

Hence, there exists a nonnegative integer $\lambda$, which satisfies
$\mathcal{Z}_k^{\lambda}(n) =1$, when $n$ is a single-digit number.

{\bf Case 2.} $n$ is a 2-digit number.

In this case, we denote $n = ck+d$, where $c \in \mathds{N}$, $d \in
\mathds{Z}_{0}^{+}$, and $c,d<k$. Thus, $\mathcal{Z}_k(n) = f(c) +
f(d)$.

Noting $f(1)=2$ and the conclusion in Case 1, under the assumption
\eqref{eq:H}, it's easy to obtain that
\begin{equation*}
\mathcal{Z}_k(n) = f(c) + f(d) \leq (c+1)+(d+1) <ck+d=n.
\end{equation*}

Therefore, there exists a nonnegative integer $\lambda$, which
satisfies $\mathcal{Z}_k^{\lambda}(n)$ is a single-digit number. It
will become Case 1.

{\bf Case 3.} $n$ is a m-digit number, where $m \geq 3$.

In this case, by Lemma 2, we can easily get that there exists a
nonnegative integer $\lambda$, $\lambda\leq m-2$, which satisfies
$\mathcal{Z}_k^{\lambda}(n)$ is a 2-digit number. It will become
Case 2, then, Case 1.

Overall, it is not difficult to see that there will exist a
nonnegative integer $\lambda$, which satisfies
$\mathcal{Z}_k^{\lambda}(n) =1$. And noting
$\mathcal{Z}_k(1)=f(1)=2$, $\mathcal{Z}_k(2)=f(2)=1$, therefore, for
any positive integer $\mu$, which satisfies $\mu \geq \lambda$, we
have $\mathcal{Z}_k^{\mu}(n) \in \{1, 2\}$.

Hence, Theorem 1 is
proved.\\
\rightline{$\Box$}

Consider a positive integer in decimalism as
\begin{equation} \label{eq:a}
\alpha=\overline{a_{m-1}a_{m-2}\dots a_{1}a_{0}},
\end{equation}
where $a_i \in \mathds{Z}_0^+$, $0 \leq a_i \leq 9$ $(0 \leq i \leq
m-1)$, and $a_{m-1} \neq 0$. As a direct result of Theorem 1, we can
easily obtain the following Corollary.

\noindent\textbf{Corollary 1}\quad \emph{For any positive integer
$\alpha$ in decimalism, that is the expression \eqref{eq:a}, the
$\mathcal{Z}$ transformation sequence
$\{\mathcal{Z}_{10}^j(\alpha),j \geq 0\}$ has a limit set $M=\{1,2\}$ when $p\in \{6,7,8\}$ in \eqref{eq:f}.}\\

\section{Examples}

\hspace{0.5cm} In this section, we give some examples for different
$n$, $k$, and $p$, to see the periodic characteristic of the
$\mathcal{Z}$ transformation.

\noindent{\bf \emph{Example 1.}} Take $n=9815671$, $k=16$, and
$p=8$, this set satisfies the assumptions \eqref{eq:H}.

\noindent{\bf \emph{Example 2.}} Take $n=71517$, $k=10$, and $p=6$,
this set satisfies the requirement of Theorem 1, also Corollary 1.

We can see from Figure 1 that by finite $\mathcal{Z}$
transformation, $n=987654321$ and $n=987654321$ will become $1$ in
their respective cases.

\noindent{\bf \emph{Example 3.}} Take $n=283$, $k=3$, and $p=2$,
this set doesn't satisfy the requirement of Theorem 1, that is
$p\ngtr 5$.

From Figure 2, it implies that $n=283$ in $3$-adic will converge to
$4$ in the sense of the $\mathcal{Z}$ transformation.

Of course, the assumptions \eqref{eq:H} in Theorem 1 are sufficient,
but not necessary. We can see this from the following example.

\noindent{\bf \emph{Example 4.}} Take $n=12345$, $k=5$, and $p=2$.
This set doesn't satisfy the requirement of Theorem 1, but we still
have that $n$ will become $1$ after six times of the $\mathcal{Z}$
transformation (see Figure 3).


\begin{figure}[htbp]
\centering \subfigure[$n=9815671$, $k=16$, $p=8$]{
\begin{minipage}[t]{0.45\linewidth}
\centering
\includegraphics[width=3.1in]{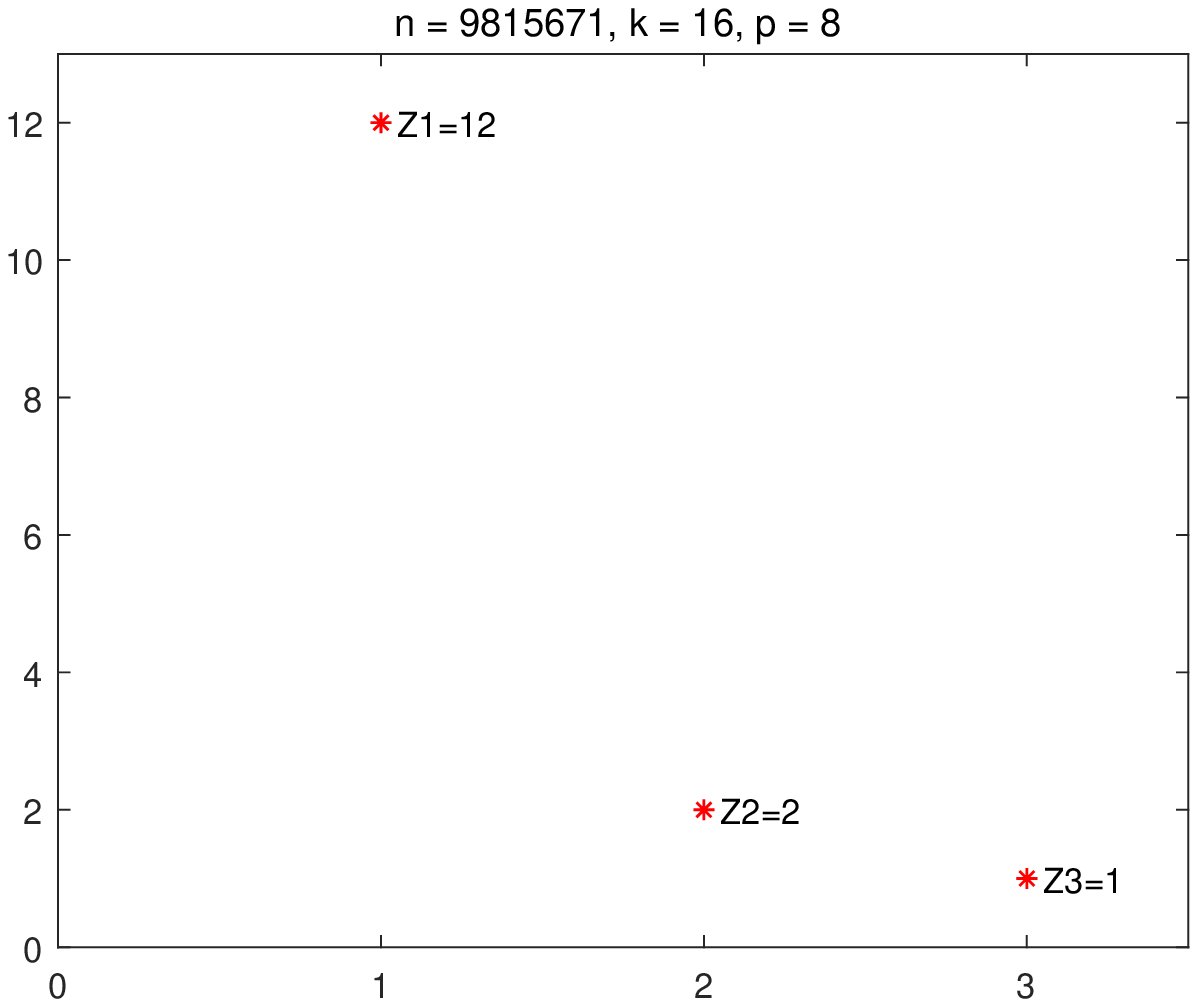}
\end{minipage}%
}%
\hspace{0.8cm} \subfigure[$n=71517$, $k=10$, $p=6$]{
\begin{minipage}[t]{0.45\linewidth}
\centering
\includegraphics[width=3.1in]{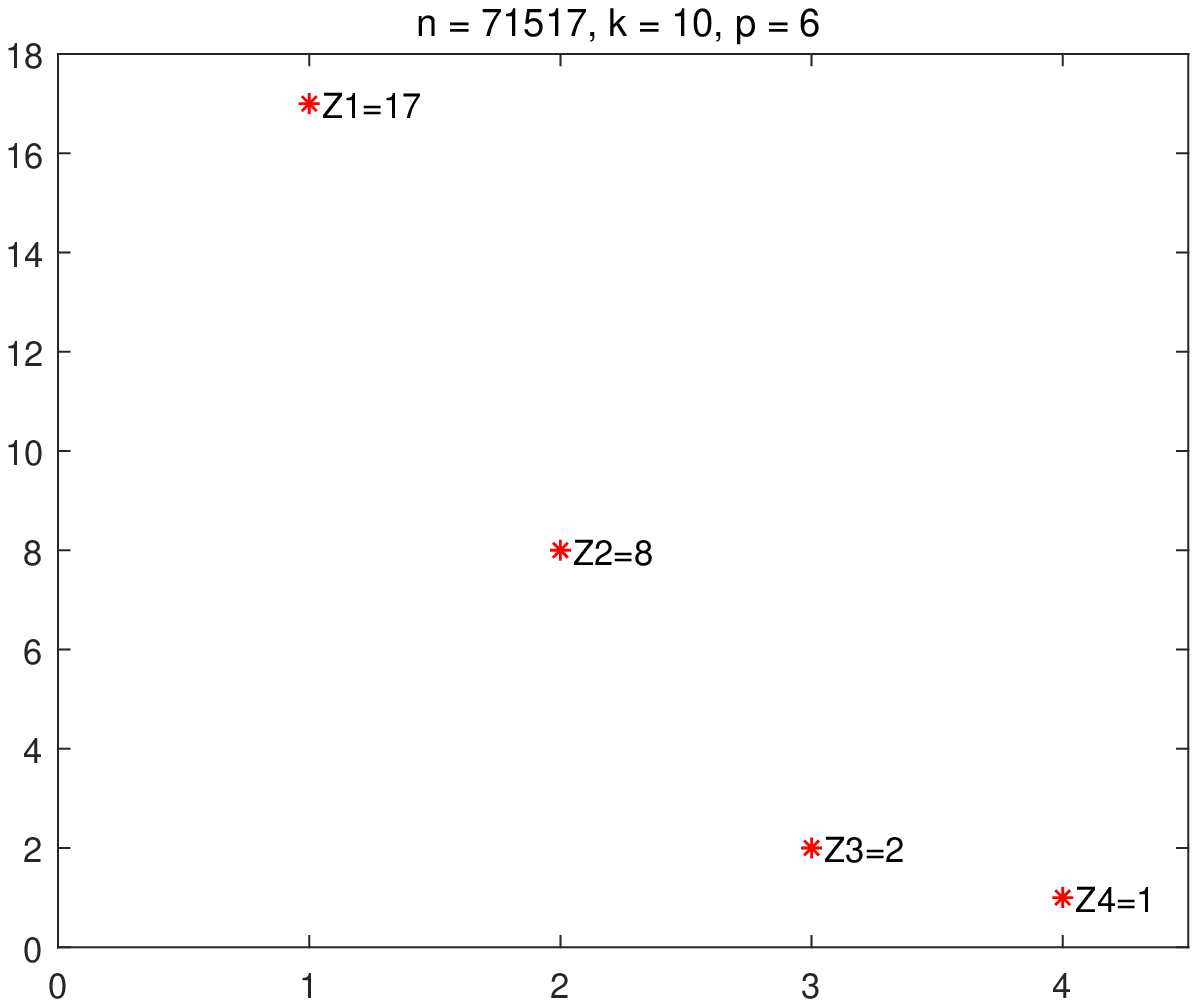}
\end{minipage}
}%
\centering \caption{Two sets of $n, k, p$, which satisfy the
assumptions \eqref{eq:H}. }
\end{figure}


\begin{figure}[htbp]
\begin{minipage}[t]{0.47\linewidth}
\includegraphics[width=3.1in, height=2.4in]{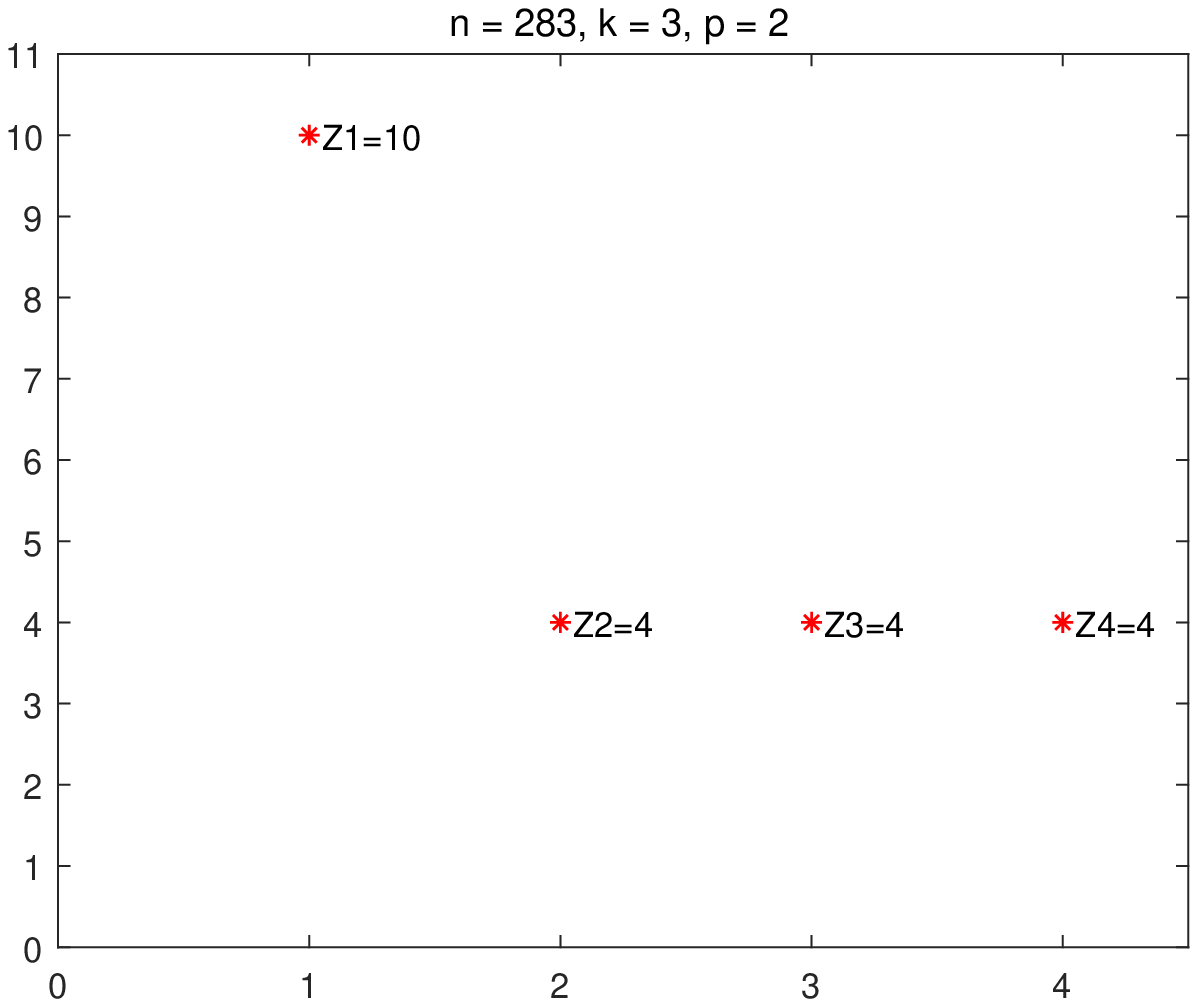}
\caption{A set of $n, k, p$, which doesn't satisfy the assumptions.}
\label{fig1}
\end{minipage}%
\hspace{0.8cm}
\begin{minipage}[t]{0.45\linewidth}
\includegraphics[width=3.1in, height=2.4in]{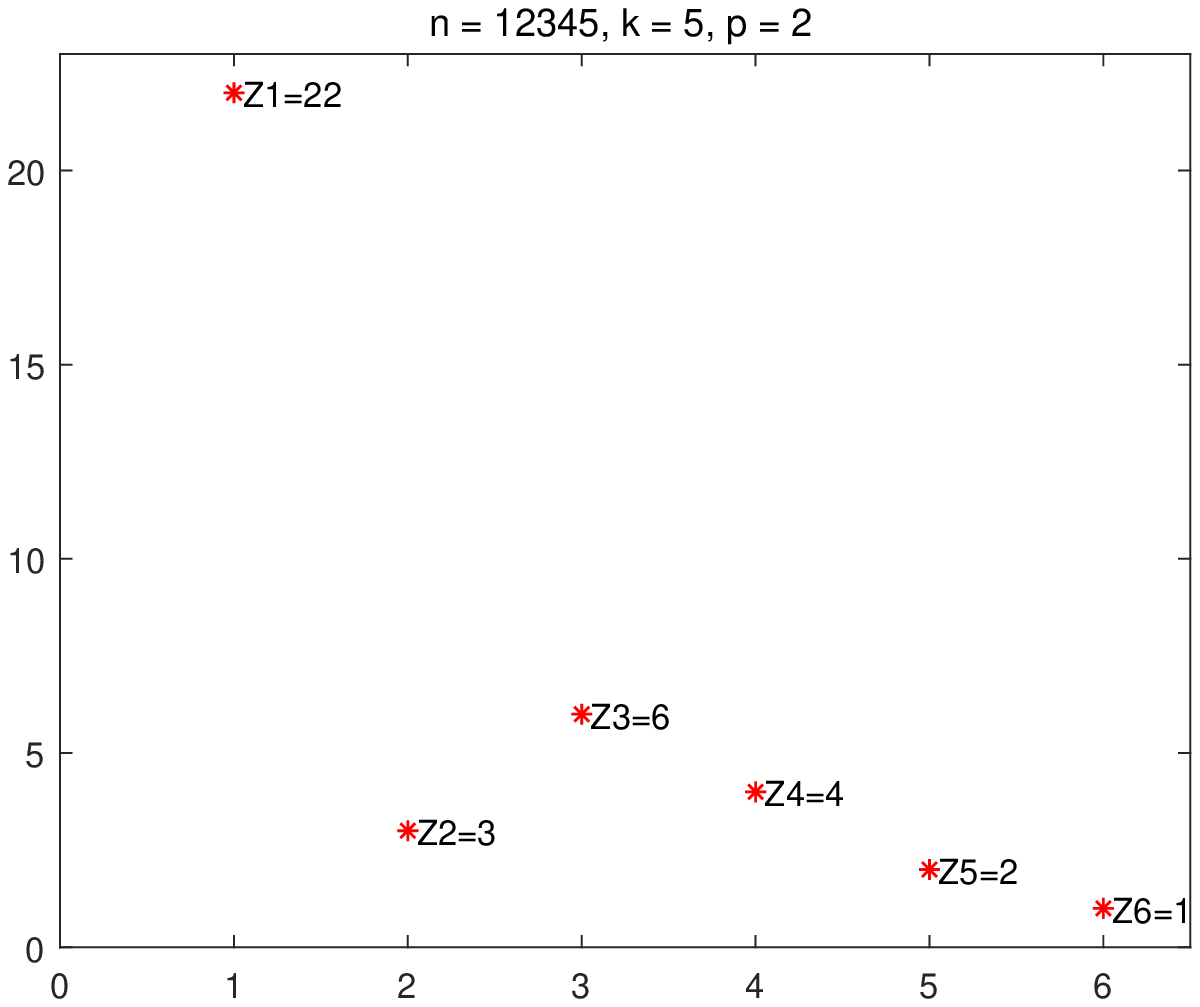}
\caption{Another set of $n, k, p$, which doesn't satisfy the
assumptions.} \label{fig2}
\end{minipage}
\end{figure}

\section*{Acknowledgements}

\hspace{0.5cm} This article is supported by Top Disciplines(Class-A)
of Zhejiang Province and Teaching Reform Project of Hangzhou Normal
University.


\end{document}